\documentclass[11pt,twoside]{article}
\usepackage{mathrsfs}
\usepackage{amsmath}
\usepackage{amssymb}
\usepackage{fancyhdr}
\usepackage{latexsym}
\usepackage{bbding}
\usepackage{mathrsfs}
\usepackage{wasysym}
\usepackage{cite}
\usepackage{xcolor}

\usepackage{multicol,graphics}

\newtheorem{theorem}{Theorem}[section]

\numberwithin{equation}{section}

\allowdisplaybreaks

\begin{document}
\title{{Notes on Liouville-type theorems for the 3D stationary Navier-Stokes equations}\thanks{\text{E-mail addresses}:
jhonglings@163.com, sunjianfeng0715@163.com, gaston.v-h@outlook.com,
 jihzhao@163.com.} }

\author{
 {\small Hongling Jiang$^{\text{1}}$,   Jianfeng Sun$^{\text{1}}$,  Gast\'on Vergara-Hermosilla$^{\text{2}}$, Jihong Zhao$^{\text{1}}$}
\\
 {\small  $^{\text{1}}$School of Mathematics and Information Science, Baoji University of Arts and Sciences,}\\
{\small   Baoji, Shaanxi 721013, China} \\
{\small $^{\text{2}}$ Institute for Theoretical Sciences, Westlake University, Hangzhou, Zhejiang 310024, China.}
}
\date{}
\maketitle

\begin{abstract}
In \cite{CV23}, Chamorro and Vergara-Hermosilla established several Liouville-type theorems to the Navier-Stokes equations in the framework of the variable Lebesgue spaces. These results may allow the variable exponent $p(\cdot)$ beyond the range of $[3,\frac{9}{2}]$ in some non-negligible regions in $\mathbb{R}^3$. In this paper we find two new non-negligible regions, in which the Liouville-type theorems still hold under some assumptions imposed on $p(\cdot)$ in these regions. Our results can be regarded as the generalization of the results in \cite{CV23}.
\end{abstract}
\smallbreak

\medskip



\section{Introduction}

\setcounter{section}{1}
\setcounter{theorem}{0}

In this paper we consider the following three dimensional incompressible stationary Navier-Stokes equations on $\mathbb{R}^3$
\begin{equation}\label{eq1.1}
		\begin{cases}
			-\Delta u+(u\cdot \nabla )u+\nabla \pi =0, \\
			\nabla \cdot u=0,
		\end{cases}
\end{equation}
where $u$ is the vector velocity field and $\pi $ is the scalar pressure.
\smallbreak

By \cite{L16}, we know the existence of solutions $(u,\pi)$ in the space $\dot{H}^1(\mathbb{R}^3) \times \dot{H}^{\frac{1}{2}}(\mathbb{R}^3)$, however, under the following additional conditions
	\begin{equation}\label{1.2}
		D(u)=\int _{\mathbb{R}^3}|\nabla u |^2dx<\infty\ \ \text{and}\ \
		\lim\limits_{|x|\to \infty } |u(x)|=0,
	\end{equation}
whether or not $u$ must be identically zero is a challenge open problem. This is commonly referred as the Liouville-type problem of equations \eqref{eq1.1}. In the monograph \cite{G11}, Galdi showed that if $u \in L^{\frac{9}{2}}(\mathbb{R}^3)$, then $u=0$. Chae and Wolf \cite{CW16} subsequently established a logarithmic improvement of Galdi's result by assuming that
	\begin{equation*}
		\int _{\mathbb{R}^3} |u|^{\frac{9}{2}}\{\text{ln}(2+|u|^{-1})\}^{-1}dx<\infty .
	\end{equation*}
Chae \cite{C14} showed the condition $\Delta u \in L^{\frac{6}{5}}(\mathbb{R}^3)$ implies $u=0$. Seregin \cite{S16} proved that if $u \in L^6(\mathbb{R}^3) \cap BMO^{-1}$, then $u=0$. Kozono, Terasawa and Wakasugi \cite{KTW17} proved that if the vorticity $\omega =\text{curl} u$ satisfies
	\begin{equation*}
		\limsup_{|x|\to \infty}|x|^{\frac{5}{3}}|\omega (x)|\le (\delta D(u))^{\frac{1}{3}}\ \ \ \text{or}\ \ \ \|u\|_{L^{\frac{9}{2},\infty}}\le (\delta ' D(u))^{\frac{1}{3}},
	\end{equation*}
then $u=0$,	where $\delta$ and $\delta '$ are two sufficiently small positive constants. Recently,   Chamorro, Jarr\'{i}n and Lemari\'{e}-Rieusset \cite{CJL21} proved that if $u \in L^p(\mathbb{R}^3)$ with $3\le p\le \frac{9}{2}$, then $u=0$. 	For further studies related to this topic we refer the readers to see \cite{C14, CNY24, J23, YX20} and the references therein.
	\smallbreak

Notice that by the classical Sobolev embeddings we have the space inclusion $\dot{H}^1(\mathbb{R}^3)\subset L^6(\mathbb{R}^3)$, which combining with the Liouville type results obtained for equations \eqref{eq1.1}, we know if $u\in\dot{H}^1(\mathbb{R}^3)\cap L^p(\mathbb{R}^3)$ with $3\leq p\leq \frac{9}{2}$, then $u=0$. Recently, Chamorro and Vergara-Hermosilla \cite{CV23} established several Liouville-type results for the equations \eqref{eq1.1} in the framework of the variable Lebesgue spaces. These results may allow the variable exponent $p(\cdot)$ beyond the range of $[3,\frac{9}{2}]$ in some non-negligible regions in $\mathbb{R}^3$, which provide us some new insights for the Liouville-type problem of equations \eqref{eq1.1}.
For additional and distinct results on Liouville theorems for the stationary Navier–Stokes equations, see \cite{V-H26a,V-H26b}.
\smallbreak

Motivated by the regions studied in \cite{CV23}, in this paper we aim to extend the Liouville theorems established there to more general non-negligible regions in $\mathbb{R}^3$.
Before we state our main results, let us introduce some notations and concepts for clarity. Let $\mathcal{P}(\Omega)$ be the set of all Lebesgue measurable functions $p(\cdot):\Omega\rightarrow[1,+\infty]$. Then for $p(\cdot)\in \mathcal{P}(\Omega)$, we define
\begin{equation}\label{Luxemburg norm}
   \|f\|_{L^{p(\cdot)}}:=\inf\Big\{\lambda>0:  \rho_{p(\cdot)}\Big(\frac{f}{\lambda}\Big)\leq 1\Big\},
\end{equation}
where the modular function $\rho_{p(\cdot)}$ associated with $p(\cdot)$ is given by the expression
\begin{equation*}
   \rho_{p(\cdot)}(f):=\int_{\Omega}|f(x)|^{p(x)}dx.
\end{equation*}
 With the Luxemburg norm \eqref{Luxemburg norm},  we define the Lebesgue space of variable exponent $L^{p(\cdot )}(\Omega)$ to be the set of Lebesgue measurable functions $f$ such that $\|f\|_{L^{p(\cdot )}}<+\infty$. Consider a measurable domain $\Omega$ such that $\Omega \subset \mathbb{R}^3$ and denote
$\mathbb{R}^3\backslash\Omega=\{x\in \mathbb{R}^3,x\notin \Omega\}$. For any $p(\cdot)\in \mathcal{P}(\mathbb{R}^3)$, we denote by $p_\Omega(\cdot )$ the variable exponent restricted to the set $\Omega$, i.e., $p_\Omega(\cdot )=p(\cdot )|_\Omega$
and denote
$$
p^{-}_{\Omega}:=\operatorname{essinf}_{x\in\Omega}p(x), \ \  p^{+}_{\Omega}:=\operatorname{esssup}_{x\in\Omega}p(x).
$$
For further background on variable Lebesgue spaces, we refer the interested reader to the books \cite{CF13,DHHR11} and/or the preliminaries sections of \cite{V-H25, V-H25b}.
\smallbreak

With this information at hand, we now state our main results. 		
\begin{theorem}\label{th1.1}
Let $G$ be a subset in $\mathbb{R}^3$ defined by
\begin{equation*}
			G=\big\{(x_1,x_2,x_3)\in \mathbb{R}^3:x_2^2+x_3^2\le x_1^{\gamma }(\ln x_1)^m,x_1>0, m\in \mathbb{Z}^+\big\},
\end{equation*}
where $0<\gamma <2$. Assume that the variable exponent $p(\cdot)\in \mathcal{P}(\mathbb{R}^3)$ satisfies the following conditions
\begin{align}\label{eq1.4}
			&3<p_{(\mathbb{R}^3\backslash G)}^-\le p_{(\mathbb{R}^3\backslash G)}(x)\le p_{(\mathbb{R}^3\backslash G)}^{+}<\frac{9}{2},\ \
			\frac{9}{2}<p_G^-\le p_{G}(x)\le p_G^+<\frac{3\gamma +3}{\gamma }.
\end{align}
For any weak solution $u\in L_{loc}^2(\mathbb{R}^3)$ and $\pi \in \mathcal{D}'(\mathbb{R}^3)$  of the stationary Navier-Stokes equations \eqref{eq1.1}, if we further assume that  $u\in L^{p(\cdot )}(\mathbb{R}^3)$  and $\pi \in L^{\frac{p(\cdot )}{2}}(\mathbb{R}^3)$, then $u=0$.
\end{theorem}
\noindent\textbf{Remark 1.1}
In \cite{CV23}, the authors considered the subset
\[
D=\left\{
(x_1,x_2,x_3)\in \mathbb{R}^3:
x_2^2+x_3^2\le x_1^{\gamma},
\ x_1>0,\
0<\gamma<1
\right\}.
\]
The set $D$ is a subset of $\mathbb{R}^3$ with infinite Lebesgue measure, and under suitable assumptions on the variable exponent $p(\cdot)$, the authors established Liouville type theorems for the stationary Navier--Stokes equations. Since
\[
D\subset G,
\]
the region considered in Theorem \ref{th1.1} is strictly more general than that in \cite{CV23}. In particular, the admissible
 region is enlarged by allowing an additional logarithmic factor in the cross-sectional growth condition. Therefore, Theorem \ref{th1.1} extends and generalizes Theorem 2 in \cite{CV23}.

\begin{theorem}\label{th1.2}
		Let $H$ be a subset in $\mathbb{R}^3$ defined by
		\begin{equation*}
			H=\big\{(x_1,x_2,x_3)\in \mathbb{R}^3: x_2^2+x_3^2\le \frac{(\ln(1+x_1))^{m}}{1+x_1}, \ x_1>0, m \in \mathbb{Z}^+\big\}.
		\end{equation*}
Assume that the variable exponent $p(\cdot)\in \mathcal{P}(\mathbb{R}^3)$ satisfies the following conditions
		\begin{align}\label{eq1.5}
			3<\mathbf{p}^-_{(\mathbb{R}^3\backslash H)}\leq \mathbf{p}_{(\mathbb{R}^3\backslash H)}(x)\leq \mathbf{p}^+_{(\mathbb{R}^3\backslash H)}<\frac{9}{2},  \ \
			\mathbf{p}_{H}(x)=+\infty.
		\end{align}
		For any weak solution $u\in L_{loc}^2(\mathbb{R}^3)$ and $\pi \in \mathcal{D}'(\mathbb{R}^3)$  of the stationary Navier-Stokes equations \eqref{eq1.1}, if we further assume that  $u\in L^{\mathbf{p}(\cdot )}(\mathbb{R}^3)$  and $\pi \in L^{\frac{\mathbf{p}(\cdot )}{2}}(\mathbb{R}^3)$, then $u=0$.
\end{theorem}
\smallbreak

\noindent\textbf{Remark 1.2}
Theorem \ref{th1.2} introduces a new class of subsets of infinite Lebesgue measure in $\mathbb{R}^3$ for which the Liouville type theorem for the stationary Navier--Stokes equations \eqref{eq1.1} remains valid. In contrast to the regions considered in \cite{CV23}, the set $H$ has a different geometric structure, since its cross-sectional width is governed by
\[
\frac{(\ln(1+x_1))^{m}}{1+x_1},
\]
which corresponds to an algebraic decay with logarithmic corrections. Therefore, Theorem \ref{th1.2} shows that the Liouville property still persists for thinner anisotropic regions involving logarithmic effects.
Moreover, the proof does not depend essentially on the particular logarithmic factor appearing in the definition of $H$. By using the same argument, the conclusion of Theorem \ref{th1.2} can be extended to more general subsets of $\mathbb{R}^3$, for example
\begin{equation*}
\widetilde{H}
=
\big\{
(x_1,x_2,x_3)\in \mathbb{R}^3:
x_2^2+x_3^2
\le
\frac{(\ln(\ln(1+x_1)))^{m}}
{(1+x_1)\ln(1+x_1)},
\ x_1>0,\
m\in\mathbb{Z}^+
\big\}.
\end{equation*}
This indicates that the Liouville type theorem remains valid for a broader family of regions whose widths shrink according to iterated logarithmic corrections.
	
	\smallbreak
	
 The proof of Theorems \ref{th1.1} and \ref{th1.2} will be given in Sections 2 and  3, respectively.

\section{Proof of Theorem \ref{th1.1}}

The main idea of the proof of Theorem \ref{th1.1} is to localize \eqref{eq1.1} by a smooth cut-off function supported in a ball $B(0,R)$, then we carefully analyze the behavior of all localized terms as $R\to +\infty$,  and in this step we will exploit the hypotheses stated in Theorem \ref{th1.1} to deduce the uniqueness of the trivial solution. For this purpose,
	let $\varphi \in C_0^{\infty }(\mathbb{R}^3)$ be a cut-off function satisfying  $0<\varphi \leq1$, with $\varphi(x)=1$ for $|x|\leq\frac{1}{2}$ and $\varphi(x)=0$ for $|x|\geq1$. For any $R>1$, we denote the rescaled function $\varphi_R(x)=\varphi (\frac{x}{R})$, thus $\varphi_R(x)=1$ if $|x|\leq \frac{R}{2}$, $\varphi_R(x)=0$ if $|x|>R$, and there exists a positive constant $C$ such that $\|\nabla \varphi _R\|_{L^\infty }\le CR^{-1}$, $\|\Delta \varphi _R\|_{L^\infty }\le CR^{-2}$.
\smallbreak

 Multiplying the first equations of  \eqref{eq1.1} by $\varphi_Ru$ and integrating over $\mathbb{R}^{3}$, by using the fact that $\operatorname{supp}(\varphi_Ru)=B_R=B(0,R)$, we get
	\begin{align}\label{NS with test function}
		\int_{B_{R}}\left[-\Delta u\cdot\left(\varphi_{R}u\right)+\left(u\cdot \nabla\right) u\cdot\left(\varphi_{R}u\right)+\nabla \pi\cdot\left(\varphi_{R}u\right)\right] dx=0.
	\end{align}
	By using the divergence free condition $\nabla \cdot u=0$, after integration by parts, all terms on the left-hand side of \eqref{NS with test function} can be rewritten as:
	\begin{align*}
		-\int_{B_{R}}\Delta u\cdot\left(\varphi_{R}u\right) dx&=\sum_{i,j=1}^{3}\int_{B_{R}}\partial _j\left (\frac{u_i^2}{2}\right )\partial _j\varphi _Rdx +\sum_{i,j=1}^{3}\int_{B_{R}}(\partial _ju_i)^2\varphi _Rdx\nonumber\\
		&=-\int_{B_{R}}\Delta \varphi_{R}\left(\frac{|u|^{2}}{2}\right)dx+\int_{B_{R}}\varphi_{R}\left|\nabla u\right|^{2}dx;
	\end{align*}
	\begin{align*}
		\int_{B_{R}}\left(u\cdot \nabla\right) u\cdot\left(\varphi_{R}u\right)dx&=\frac{1}{2}\sum_{i,j=1}^{3}\int_{B_{R}}\partial_{i}\left(u_{i}u_{j}^{2}\right)\varphi_{R}dx=-\frac{1}{2}\int_{B_{R}}\nabla\varphi_{R}\cdot\left(|u|^{2}u\right)dx
	\end{align*}
	and
	\begin{align*}
		\int_{B_{R}}\nabla \pi\cdot\left(\varphi_{R}u\right)dx=-\int_{B_{R}}\nabla\varphi_{R}\cdot(\pi u)dx.
	\end{align*}
	Taking all above identities into \eqref{NS with test function}, we obtain
	\begin{align}\label{eq3.2}
		\int_{B_{R}}\varphi_{R}\left|\nabla u\right|^{2}dx &=\int_{B_{R}}\Delta \varphi_{R}\left(\frac{|u|^{2}}{2}\right)+\nabla\varphi_{R}\cdot\left(\frac{|u|^{2}}{2}u\right)
		+\nabla\varphi_{R}\cdot(\pi u)dx:=I_1+I_2+I_3.
	\end{align}
Since $\varphi _R(x)=1$ over the set $|x|\leq\frac{R}{2}$, one can easily see that
	\begin{equation*}
		\int_{B_{\frac{R}{2}}}{|\nabla u|^2}dx\le \sum\limits_{i=1}^3{|I_i|}.
	\end{equation*}
	Thus if we can prove all the limits
	\begin{equation}\label{eq3.6}
		\lim\limits_{R\rightarrow+\infty}|I_i|=0\,(i=1,2,3),
	\end{equation}
	then we obtain
	\begin{equation*}
		\lim\limits_{R\rightarrow+\infty}\int_{B_{\frac{R}{2}}}{|\nabla u|^2}dx=\|u\|_{\dot{H}^1}^2=0,
	\end{equation*}
	which implies that $\|u\|_{L^{6}}=0$ by the Sobolev embedding $\dot{H}^1(\mathbb{R}^3)\hookrightarrow L^6(\mathbb{R}^3)$, and we get $u=0$. 	
	Now we estimate the terms $I_i$ ($i=1,2,3$) one by one to complete the proof of \eqref{eq3.6}.
\smallbreak
	
	\textit{Estimate of $I_1$.} Applying the H\"older's inequality with $\frac{2}{p(\cdot )}+\frac{1}{q(\cdot )}=1$, one sees that
	\begin{align}\label{eq2.7}
		I_{1}\le \int_{\mathbb{R}^3}\left|\Delta \varphi _R\right|\frac{|u|^2}{2}dx\le C\|\Delta \varphi _R\|_{L^{q(\cdot )}(\mathbb{R}^3)}\|u\|^2_{L^{p(\cdot )}(\mathbb{R}^3)}.
	\end{align}
	Under the assumptions of Theorem \ref{th1.1}, we know that $\|u\|_{L^{p(\cdot )}(\mathbb{R}^3)}<+\infty$,  thus it suffices to estimate the quantity $\|\Delta \varphi _R\|_{L^{q(\cdot )}(\mathbb{R}^3)}$ in \eqref{eq2.7}. Since the variable exponent $p(\cdot)$ satisfies the condition \eqref{eq1.4}, it is easy to verify that the variable exponent $q(\cdot )=\frac{p(\cdot )}{p(\cdot )-2}$ satisfies
\begin{equation}\label{eq2.08}
		\frac{9}{5}<q_{(\mathbb{R}^3\backslash G)}^{-}\le q_{(\mathbb{R}^3\backslash G)}(x)\le q_{(\mathbb{R}^3\backslash G)}^{+}<3,\ \
		\frac{3\gamma+3}{\gamma+3}<q_G^- \le q_{G}(x)\le q_G^+<\frac{9}{5}.
	\end{equation}
  Since  $\operatorname{supp}(\Delta \varphi _R)\subset C(\frac{R}{2},R)=\{x\in{R^3}:\frac{R}{2}\le |x|\le R\}$, we have
	\begin{equation*}
		\|\Delta \varphi _R\|_{L^{q(\cdot )}(\mathbb{R}^3)}=\|\Delta \varphi _R\|_{L^{q(\cdot )}(C(\frac{R}{2},R))}.
	\end{equation*}
To continue, we denote by $G_1$ and $G_2$ the subsets of $G$ as $G_1:=C(\frac{R}{2},R)\cap G$ and $G_2:=C(\frac{R}{2},R)\backslash G$,
thus
	\begin{align}\label{eq2.10}
		\|\Delta \varphi _R\|_{L^{q(\cdot )}(C(\frac{R}{2},R))}
		&=\|\Delta {\varphi _R}(1_{G_1}+1_{G_2})\|_{L^{q(\cdot )}(C(\frac{R}{2},R))}\nonumber\\
    &\le\| \Delta \varphi _R\|_{L^{q(\cdot )}(G_1)}+\|\Delta \varphi _R\|_{L^{q(\cdot )}(G_2)}.
	\end{align}
For the first term on the right-hand side of \eqref{eq2.10}, we can bound it as
	\begin{align}\label{eq2.11}
		\|\Delta \varphi_R\|_{L^{q(\cdot )}(G_1)}&\le \|\Delta \varphi_R\|_{L^{\infty }(G_1)}\|1\|_{L^{q(\cdot )}(G_1)}\nonumber\\
    &\le \|\Delta \varphi_R\|_{L^{\infty }(C(\frac{R}{2},R))}\|1\|_{L^{q(\cdot )}(G_1)}\nonumber\\
    &\le C{R^{-2}}\|1\|_{L^{q(\cdot)}(G_1)},
	\end{align}
	where we have used the facts $G_1\subset C(\frac{R}{2},R)$ and $\|\Delta \varphi _R\|_{L^\infty }\le CR^{-2}$.
As $G_1\subset G$, we know $q^-_G\le q^-_{G_1}\le q^+_{G_1}\le q^+_{G}$, thus applying the Lemma 2.1 in \cite{CV23} leads to
	\begin{align}\label{h4.1}
		\|\Delta \varphi_R\|_{L^{q(\cdot)}(G_1)}
		\leq CR^{-2} \max \Big\{ |G_1|^{\frac{1}{q_G^-}}, |G_1|^{\frac{1}{q_G^+}} \Big\}\leq CR^{-2}|G_1|^{\frac{1}{q_G^-}}.
	\end{align}
On the other hand, notice that
	\begin{equation*}
		G_1\subset \mathcal{A}=\left\{(x_1,x_2,x_3)\in \mathbb{R}^3:x_2^2+x_3^2\le x_1^{\gamma }(\ln x_1)^m, 0<x_1<R\right\},
	\end{equation*}
	which yields that the volume of $G_1$ can be controlled by the volume of $\mathcal{A}$, i.e.,
	\begin{align}\label{E}
		|G_1|&\le |\mathcal{A}|= C\int_{0}^{R}x_1^\gamma (\ln x_1)^mdx_1=C\int _{0}^{R}(\ln x_1)^md\Big(\frac{x_1^{\gamma+1}}{\gamma+1}\Big)\nonumber\\
		&=C\Big(\Big.{\frac{x_1^{\gamma+1}(\ln x_1)^m}{\gamma+1} }\Big|_0^R -\frac{m}{\gamma+1}\int_0^Rx_1^\gamma (\ln x_1)^{m-1}dx_1 \Big)\nonumber\\
		&=C\Big(\frac{R^{\gamma+1}(\ln R)^m}{\gamma+1}-\frac{mR^{\gamma+1}(\ln R)^{m-1}}{(\gamma+1)^2}+\cdots+(-1)^{m+1}\frac{m!R^{\gamma+1}}{(\gamma+1)^{m+1}}\Big).
	\end{align}
Taking \eqref{E} into \eqref{h4.1} yields
	\begin{equation}\label{h4.3}
		\|\Delta \varphi_R\|_{L^{q(\cdot )}(G_1)}\le C{R^{-2}}\Big(\frac{R^{\gamma+1}(\ln R)^m}{\gamma+1}-\frac{mR^{\gamma+1}(\ln R)^{m-1}}{(\gamma+1)^2}+\cdots+(-1)^{m+1}\frac{m!R^{\gamma+1}}{(\gamma+1)^{m+1}}\Big)^{\frac{1}{q_G^-}}.
	\end{equation}
	Since $m\in \mathbb{Z}^+$, it suffices to consider only the first term on the right-hand side of \eqref{h4.3},
	\begin{equation*}
		\lim\limits_{R\rightarrow +\infty}C{R^{-2}}\Big(\frac{R^{\gamma+1}(\ln R)^m}{\gamma+1}\Big)^{\frac{1}{q_G^-}}\le\lim\limits_{R\rightarrow +\infty} CR^{-2+\frac{\gamma +1}{q_G^-}}(\ln R)^{\frac{m}{q_G^-}}.
	\end{equation*}
By \eqref{eq2.08} and $0 < \gamma < 2$, we know $ \frac{1}{q_G^-} < \frac{\gamma+3}{3(\gamma+1)}$, thus $-2 + \frac{\gamma+1}{q_G^-} < 0$, and we obtain
	\begin{equation}\label{eq2.14}
		\lim\limits_{R \to +\infty} \|\Delta \varphi_R\|_{L^{q(\cdot)}(G_1)} = 0,
	\end{equation}
	where we used the fact that $\lim\limits_{R \to +\infty}R^\alpha (\ln R)^m=0$ for $\alpha<0$ and any $m\in \mathbb{Z}^+$.
\smallbreak

For the second term on the right-hand side of \eqref{eq2.10},  applying  the same arguments,  one sees that
	\begin{align}\label{eq2.15}
		\|\Delta \varphi_R\|_{L^{q(\cdot )}(G_2)}&\le \|\Delta \varphi_R\|_{L^{\infty }(G_2)}\|1\|_{L^{q(\cdot )}(G_2)}\nonumber\\
    &\le C{R^{-2}}\|1\|_{L^{q(\cdot)}(G_2)}\nonumber\\
    &\le CR^{-2}\max \Big\{|G_2|^{\frac{1}{q_{G_2}^-}},|G_2|^{\frac{1}{q_{G_2}^+}}\Big\}.
	\end{align}
By \eqref{eq2.08} again, we know that
$\frac{1}{3}<\frac{1}{q_{\mathbb{R}^3\backslash G}^{+}(x)}<\frac{1}{q_{G_2}^{+}(x)}\leq \frac{1}{q_{G_2}^{-}(x)}\leq \frac{1}{q_{\mathbb{R}^3\backslash G}^{-}(x)}<\frac{5}{9}$, and $|G_2| \leq |C(\frac{R}{2}, R)| \leq CR^3$,   which yields that
	\begin{equation}\label{eq2.16}
		\|\Delta \varphi_R\|_{L^{q(\cdot)}(G_2)}\le C\max \Big\{R^{-2+\frac{3}{q_{\mathbb{R}^3\backslash G}^-}}, \  R^{-2+\frac{3}{q_{\mathbb{R}^3\backslash G}^+}}\Big\}.
	\end{equation}
Since $-2+\frac{3}{q_{\mathbb{R}^3\backslash G}^{+}} \leq -2+\frac{3}{q_{\mathbb{R}^3\backslash G}^{-}} <0$,  we obtain
	\begin{equation}\label{eq2.17}
		\lim\limits_{R\rightarrow +\infty}\|\Delta \varphi_R\|_{L^{q(\cdot )}(G_2)}= 0.
	\end{equation}
	Putting \eqref{eq2.14} and \eqref{eq2.17} together,   we conclude from \eqref{eq2.10} that $\lim\limits_{R\rightarrow +\infty}\|\Delta \varphi _R\|_{L^{q(\cdot )}(C(\frac{R}{2},R))}=0.$
	Now by the assumption $\|u\|_{L^{p(\cdot )}(\mathbb{R}^3)}<+\infty,$ we obtain from \eqref{eq2.7} that $\lim\limits_{R\rightarrow +\infty}|I_1|=0$.
\smallbreak
	
	\textit{Estimate of $I_2$.} Using the H\"older's inequality  yields that
	\begin{align}\label{eq2.21}
		|I_2|\le \frac{1}{2}\int_{C(\frac{R}{2},R)}{|u|^3|\nabla \varphi _R|}dx
		\le C\|\nabla \varphi_R\|_{L^{r(\cdot )}(C(\frac{R}{2},R))}\|u\|^3_{L^{p(\cdot )}(\mathbb{R}^3)},
	\end{align}
	where $\frac{3}{p(\cdot )}+\frac{1}{r(\cdot )}=1$.  By \eqref{eq1.4}, we deduce that the variable exponent $r(\cdot )=\frac{p(\cdot )}{p(\cdot )-3}$ satisfies
	\begin{equation}\label{eq2.22}
		3<r_{(\mathbb{R}^3\backslash G)}^{-}\le r_{(\mathbb{R}^3\backslash G)}(x)\le r_{(\mathbb{R}^3\backslash G)}^{+}<+\infty,\ \
		\gamma +1<r_{G}^{-}\le r_{G}(x)\le r_{G}^{+}<3.
	\end{equation}
By the definition of the sets $G_1$ and $G_2$ and proceeding just as  \eqref{eq2.10}, one has
	\begin{align}\label{eq3.55555}
		\|\nabla \varphi _R\|_{L^{r(\cdot )}(C(\frac{R}{2},R))} \leq \| \nabla \varphi _R\|_{L^{r(\cdot )}(G_1)} + \|\nabla  \varphi _R\|_{L^{r(\cdot )}(G_2)}.
	\end{align}
		Following the same ideas  that leads to estimates \eqref{h4.3} and \eqref{eq2.16} (with the difference that $\|\nabla \varphi_R\|_{L^{\infty }} \leq CR^{-1}$), we obtain that
	\begin{align}\label{eq2.23}
		\|\nabla \varphi_R\|_{L^{r(\cdot )}(C(\frac{R}{2},R))}
		&\le C{R^{-1}}\Big(\frac{R^{\gamma+1}(\ln R)^m}{\gamma+1}-\frac{mR^{\gamma+1}(\ln R)^{m-1}}{(\gamma+1)^2}+\cdots+(-1)^{m+1}\frac{m!R^{\gamma+1}}{(\gamma+1)^{m+1}}\Big)^{\frac{1}{r_G^-}}\nonumber\\
		&+ C\max \Big\{R^{-1+\frac{3}{r_{\mathbb{R}^3\backslash G}^{-}}}, R^{-1+\frac{3}{r_{\mathbb{R}^3\backslash G}^+}}\Big\}.
	\end{align}
Notice that
	\begin{equation*}
		\lim\limits_{R\rightarrow +\infty}C{R^{-1}}\Big(\frac{R^{\gamma+1}(\ln R)^m}{\gamma+1}\Big)^{\frac{1}{r_G^-}}\le\lim\limits_{R\rightarrow +\infty} CR^{-1+\frac{\gamma +1}{r_G^-}}(\ln R)^{\frac{m}{r_G^-}}.
	\end{equation*}
Moreover, by the range of the exponent $r(\cdot )$ in \eqref{eq2.22}, one sees that $-1+\frac{\gamma +1}{r_G^-}<0$  and $-1+\frac{3}{r_{\mathbb{R}^3\backslash G}^{+}}\le -1+\frac{3}{r_{\mathbb{R}^3\backslash G}^{-}}<0$.
	This readily implies that
	\begin{equation}\label{eqg2}
		\lim\limits_{R\rightarrow +\infty}\|\nabla \varphi_R\|_{L^{r(\cdot )}(C(\frac{R}{2},R))}=0.
	\end{equation}
	Since $\|u\|_{L^{p(\cdot )}(\mathbb{R}^3)}<+\infty$, we obtain from \eqref{eq2.21} that $\lim\limits_{R\rightarrow +\infty}|I_2|=0$.
\smallbreak
	
	\textit{Estimate of $I_3$}. Since $u\in L^{p(\cdot )}(\mathbb{R}^3)$ and $\pi \in L^{\frac{p(\cdot )}{2}}(\mathbb{R}^3)$, we can get the limit $\lim\limits_{R\rightarrow +\infty}|I_3|=0$ by \eqref{eqg2} and the following inequality ($\frac{3}{p(\cdot )}+\frac{1}{r(\cdot )}=1$):
	\begin{align}\label{eq2.24}
		|I_3|&\le \int_{C(\frac{R}{2},R)}{|\pi||u||\nabla \varphi_R|}dx\le C\|\nabla \varphi_R\|_{L^{r(\cdot )}(C(\frac{R}{2},R))}\|\pi \|_{L^{\frac{p(\cdot )}{2}}(\mathbb{R}^3)}\|u\|_{L^{p(\cdot )}(\mathbb{R}^3))}.
	\end{align}
\smallbreak

		We have thus proven that, under the assumptions of Theorem \ref{th1.1}, all the terms $|I_i|$ $(i=1,2,3)$ given in \eqref{eq3.2} tend to $0$ as $R\rightarrow +\infty $, the proof of  Theorem \ref{th1.1} is completed.

\section{Proof of Theorem \ref{th1.2}}

Following the main ideas used in the proof of Theorem \ref{th1.1}, it suffices to prove that, under the assumptions of Theorem \ref{th1.2}, all terms $|I_i|$ $(i=1,2,3)$ given in \eqref{eq3.6} tend to $0$ as $R\rightarrow +\infty $.
	\smallbreak

	\textit{Estimate of $I_1$.} Applying the H\"older's inequality with $\frac{2}{\mathbf{p}(\cdot )}+\frac{1}{\mathbf{q}(\cdot )}=1$, one sees that
	\begin{align}\label{eq3.7}
		I_{1}\le C\|\Delta \varphi _R\|_{L^{\mathbf{q}(\cdot )}(C(\frac{R}{2},R))}\|u\|^2_{L^{\mathbf{p}(\cdot )}(\mathbb{R}^3)},
	\end{align}
 By \eqref{eq1.5}, we know the variable exponent $\mathbf{q}(\cdot )=\frac{\mathbf{p}(\cdot )}{\mathbf{p}(\cdot )-2}$ satisfies
	\begin{equation}\label{eq3.8}
		\frac{9}{5}<\mathbf{q}_{(\mathbb{R}^3\backslash H)}^{-}\le \mathbf{q}_{(\mathbb{R}^3\backslash H)}(x)\le \mathbf{q}_{(\mathbb{R}^3\backslash H)}^{+}<3,\ \
		\mathbf{q}_H^- =\mathbf{q}_{H}(x)=\mathbf{q}_H^+=1.
	\end{equation}
We denote by $H_1$ and $H_2$ the subsets of $H$ as  $H_1:=C(\frac{R}{2},R)\cap H$ and $H_2:=C(\frac{R}{2},R)\backslash H$,
thus
	\begin{align}\label{eq3.10}
		\|\Delta \varphi _R\|_{L^{\mathbf{q}(\cdot )}(C(\frac{R}{2},R))}
		\le\| \Delta \varphi _R\|_{L^{\mathbf{q}(\cdot )}(H_1)}+\|\Delta \varphi _R\|_{L^{\mathbf{q}(\cdot )}(H_2)}.
	\end{align}
	For the first term on the right-hand side of \eqref{eq3.10}, we can bound it as
	\begin{align}\label{eq3.11}
		\|\Delta \varphi_R\|_{L^{\mathbf{q}(\cdot )}(H_1)}&\le \|\Delta \varphi_R\|_{L^{\infty }(H_1)}\|1\|_{L^{\mathbf{q}(\cdot )}(H_1)}\nonumber\\
    &\le C{R^{-2}}\|1\|_{L^{\mathbf{q}(\cdot)}(H_1)}\nonumber\\
    &\le C{R^{-2}}\max  \big\{|H_1|^{\frac{1}{\mathbf{q}_H^-}},|H_1|^{\frac{1}{\mathbf{q}_H^+}}\big\}\nonumber\\
		&\le C{R^{-2}}|H_1|,
	\end{align}
	 where we have used  $\mathbf{q}_H^-= \mathbf{q}_H^+=1$. Notice that
	\begin{equation*}
		H_1\subset \mathcal{B}=\big\{(x_1,x_2,x_3)\in \mathbb{R}^3:x_2^2+x_3^2\le \frac{1}{1+x_1}(\ln(1+x_1))^{m},0<x_1<R\big\},
	\end{equation*}
	which yields that the volume of $H_1$ can be controlled by the volume of $\mathcal{B}$, i.e.,
	\begin{equation}\label{HH}
		|H_1|\le |\mathcal{B}|\le C\int_{0}^{R}{\frac{1}{1+x_1}(\ln(1+x_1))^{m}}dx_1=\frac{C}{m+1}(\ln(1+R))^{m+1}.
	\end{equation}
By \eqref{eq3.11} and \eqref{HH}, we have
	\begin{equation}\label{eq3.13}
		\|\Delta \varphi_R\|_{L^{\mathbf{q}(\cdot )}(H_1)}\le C{R^{-2}}(\ln(1+R))^{m+1},
	\end{equation}
	and for any $m\in \mathbb{Z}^+$, we obtain
	\begin{equation}\label{eq3.14}
		\lim\limits_{R\rightarrow +\infty}\|\Delta \varphi_R\|_{L^{\mathbf{q}(\cdot)}(H_1)}=0.
	\end{equation}
	For the second term on the right-hand side of \eqref{eq3.10},  applying  the same arguments,  one sees that
	\begin{equation}\label{eq3.16}
		\|\Delta \varphi_R\|_{L^{\mathbf{q}(\cdot)}(H_2)}\le C\max \big\{R^{-2+\frac{3}{\mathbf{q}_{\mathbb{R}^3\backslash H}^-}}, \  R^{-2+\frac{3}{\mathbf{q}_{\mathbb{R}^3\backslash H}^+}}\big\}.
	\end{equation}
	Observing that $-2+\frac{3}{\mathbf{q}_{\mathbb{R}^3\backslash H}^{+}} \leq -2+\frac{3}{\mathbf{q}_{\mathbb{R}^3\backslash H}^{-}} < 0$,  which yields that
	\begin{equation}\label{eq3.17}
		\lim\limits_{R\rightarrow +\infty}\|\Delta \varphi_R\|_{L^{\mathbf{q}(\cdot )}(H_2)}= 0.
	\end{equation}
	Taking \eqref{eq3.14} and \eqref{eq3.17} back to \eqref{eq3.10}, we obtain $ \lim\limits_{R\rightarrow +\infty}\|\Delta \varphi _R\|_{L^{\mathbf{q}(\cdot )}(C(\frac{R}{2},R))}=0,$
	and since $\|u\|_{L^{\mathbf{p}(\cdot )}(\mathbb{R}^3)}<+\infty$, we obtain from \eqref{eq3.7} that $\lim\limits_{R\rightarrow +\infty}|I_1|=0$.
\smallbreak
	
	\textit{Estimate of $I_2$.} Using the H\"older's inequality with  $\frac{3}{\mathbf{p}(\cdot )}+\frac{1}{\mathbf{r}(\cdot )}=1$, one has
	\begin{align}\label{eq3.21}
		|I_2|\le C\|\nabla \varphi_R\|_{L^{\mathbf{r}(\cdot )}(C(\frac{R}{2},R))}\|u\|^3_{L^{\mathbf{p}(\cdot )}(\mathbb{R}^3)},
	\end{align}
By \eqref{eq1.5}, we deduce that the variable exponent $\mathbf{r}(\cdot )=\frac{\mathbf{p}(\cdot )}{\mathbf{p}(\cdot )-3}$ satisfies
	\begin{equation}\label{eq3.22}
		3<\mathbf{r}_{(\mathbb{R}^3\backslash H)}^{-}\le \mathbf{r}_{(\mathbb{R}^3\backslash H)}(x)\le \mathbf{r}_{(\mathbb{R}^3\backslash H)}^{+}<+\infty \ \ \mathrm{and}\ \
		\mathbf{r}_{H}^{-}= \mathbf{r}_{H}(x)= \mathbf{r}_{H}^{+}=1,
	\end{equation}
thus we can follow the same ideas used in the estimates of \eqref{eq3.13} and \eqref{eq3.16} to obtain that
	\begin{align*}
		\|\nabla \varphi_R\|_{L^{\mathbf{r}(\cdot )}(C(\frac{R}{2},R))}&\le \|\nabla \varphi_R\|_{L^{\mathbf{r}(\cdot )}(H_1)}+\|\nabla \varphi _R\|_{L^{\mathbf{r}(\cdot )}(H_2)}\nonumber\\
		&\le C{R^{-1}}(\ln(1+R))^{1+m}+C\max \big\{R^{-1+\frac{3}{\mathbf{r}_{\mathbb{R}^3\backslash H}^{-}}}, R^{-1+\frac{3}{\mathbf{r}_{\mathbb{R}^3\backslash H}^+}}\big\},
	\end{align*}
 which by the range of $\mathbf{r}(\cdot )$ in \eqref{eq3.22}, one sees that $-1+\frac{3}{\mathbf{r}_{\mathbb{R}^3\backslash H}^{+}}\le -1+\frac{3}{\mathbf{r}_{\mathbb{R}^3\backslash H}^{-}}<0.$	This readily implies that
	\begin{equation}\label{eq3.25}
		\lim\limits_{R\rightarrow +\infty}\|\nabla \varphi_R\|_{L^{\mathbf{r}(\cdot )}(C(\frac{R}{2},R))}=0.
	\end{equation}
By the assumption $u \in L^{\mathbf{p}(\cdot )}(\mathbb{R}^3)$, we know from \eqref{eq3.21} that $\lim\limits_{R\rightarrow +\infty}|I_2|=0$.
\smallbreak
	
	\textit{Estimate of $I_3$.} Since $u\in L^{\mathbf{p}(\cdot )}(\mathbb{R}^3)$ and $\pi \in L^{\frac{\mathbf{p}(\cdot )}{2}}(\mathbb{R}^3)$, we can get the limit $\lim\limits_{R\rightarrow +\infty}|I_3|=0$ by \eqref{eq3.25} and the following inequality ($\frac{3}{\mathbf{p}(\cdot )}+\frac{1}{\mathbf{r}(\cdot )}=1$):
	\begin{align*}\label{eq3.24}
		|I_3|&\le \int_{C(\frac{R}{2},R)}{|\pi||u||\nabla \varphi_R|}dx\le C\|\nabla \varphi_R\|_{L^{\mathbf{r}(\cdot )}(C(\frac{R}{2},R))}\|\pi \|_{L^{\frac{\mathbf{p}(\cdot )}{2}}(\mathbb{R}^3)}\|u\|_{L^{\mathbf{p}(\cdot )}(\mathbb{R}^3))}.
	\end{align*}
	\smallbreak

	We have thus proven that, under the assumptions of Theorem \ref{th1.2}, all the terms $I_{i} (i=1,2,3)$ given in \eqref{eq3.6} to $0$ as $R \to +\infty$. The proof of Theorem \ref{th1.2} is completed.
\bigskip
		
		\noindent \textbf{Acknowledgements}
 The authors declared that they have no conflict of interest. This work is partially supported by the National Natural Science Foundation of China (no. 12361034) and the Graduate Student Innovation Project of Baoji University of Arts and Sciences (no. YJSCX25YB38).
		\vskip.2cm

	\noindent \textbf{Data Availability Statement} No data was used for the research described in the article.

\small{

}

\end{document}